\documentclass[twocolumn]{IEEEtran}

\usepackage{amsthm,amsmath, amssymb, graphicx, bbm, mathrsfs, color}

\newtheorem{lemma}{Lemma}

\theoremstyle{definition}

\newtheorem{remark}{Remark}
\newtheorem{example}{Example}

\title{Differentiation via Logarithmic Expansions}

\author{Michael C. Fu$^\dagger$, Bernd Heidergott$^\sharp$, Haralambie Leahu$^\$$ and Felisa V{\'a}zquez-Abad$^\ddag$
\thanks{
$^\dagger$ Smith School of Business \& Institute for Systems Research, University of Maryland, College Park,
E-Mail: mfu@umd.edu,

$\sharp$ Department of Econometrics and Operations Research, VU Amsterdam,
E-Mail: bheidergott@feweb.vu.nl, 

$\$$ Department of Department of Mathematics, University of Amsterdam,
E-Mail: haralambie@gmail.com,

$\ddag$ Hunter, College, City University of New York,
 E-Mail: felisav@hunter.cuny.edu.
}
}

\begin{document}

\maketitle

\begin{abstract}
In this note, we introduce a new finite difference approximation 
called the Black-Box Logarithmic Expansion Numerical Derivative (BLEND) algorithm, 
which is based on a formal logarithmic expansion of the differentiation operator. 
BLEND capitalizes on parallelization and provides derivative approximations of arbitrarily precision, 
i.e., our analysis can be used to determine the number of terms in the series expansion to guarantee a specified number of decimal places of accuracy.  
Furthermore, in the vector setting, the complexity of the resulting directional derivative
 is independent of the dimension of the parameter. 
\\

\end{abstract}

\section{Introduction}

Evaluation of derivatives is essential in optimization, control and sensitivity analysis. 
We consider the setting where the function of interest, 
$\phi:\Theta\rightarrow\mathbb{R},~ \Theta \subset \mathbb{R} $,
is not available in closed-form but function evaluations are available for any $\theta \in \Theta$,  
i.e., we are in a ``black box" scenario where derivatives of $ \phi $ have to be computed numerically by means of a finite difference (FD) approximation.  
For instance, the stationary distribution of a finite Markov chain may fail to have a closed-form solution,
which rules out analytical computation of derivatives, 
but the stationary distribution can be easily evaluated numerically. 
Many complex queueing networks fall into this category, and we use a simple tandem queue to illustrate our algorithm.

The most well-known FD approximations, with $ h > 0$, are the forward and central FD approximations given respectively by
\begin{equation}\label{eq:forwardFD}
\phi'(\theta)  \approx \frac{ \phi ( \theta + h ) - \phi  ( \theta )}{ h } ,
\end{equation}
and  
\begin{equation}\label{eq:centralFD}
\phi'(\theta)  \approx \frac{ \phi ( \theta + h ) - \phi  ( \theta  - h )}{2 h} .
\end{equation}
Many other FD approximations involving weighted sums and differences of function evaluations on one-dimensional grids with arbitrary spacing
are available in the literature, including for higher-order derivatives; see, for example, \cite{1}.
However, determining the weights required for each grid  point is often computationally demanding and depends on the chosen grid. 
Also, most of the approximations involve recomputing the weights when additional precision is required.  

In this note, we consider a simple FD series approximation that we call  the Black-Box Logarithmic Expansion Numerical Derivative (BLEND), which is based on a formal logarithmic expansion of the differentiation operator. 
BLEND provides derivative approximations of arbitrary precision, in the sense that  a specified number of decimal places of accuracy can be guaranteed by adding a sufficient number of terms in the series approximation, so the stopping criterion can be evaluated directly as the algorithm runs.  
The proposed BLEND algorithm capitalizes on modern computing platforms 
by exploiting parallel processing in evaluating the terms simultaneously.

BLEND begins with the Taylor series expansion:
\begin{equation}
\phi(\theta+h)=\sum_{n=0}^\infty\frac{h^n}{n!}\:\phi^{(n)}(\theta), ~~\forall h>0.
\label{eq.Taylor}
\end{equation}
Denoting by $\Delta$ the differential operator $\Delta\phi:=\phi'$ and by $\mathcal{T}_h$ the shift-operator $\mathcal{T}_h\phi:=\phi(\cdot+h)$,
(\ref{eq.Taylor}) can be rewritten in operator notation as 
$$
\mathcal{T}_h\phi =\sum_{n=0}^\infty\frac{h^n}{n!}\:\Delta^n\phi 
 = \sum_{n=0}^\infty\frac{(h\Delta)^n\phi}{n!} .
$$
leading to the compact operator relationship
$\mathcal{T}_h=\exp(h\Delta)$.
The existence of series expansions for inverses of analytic functions is a well-known result
of complex analysis (see \cite{nnn,mmm}), and formal inversion yields
\begin{equation}
\Delta=\frac{\ln(\mathcal{T}_h)}{h},
%\Delta=(1/h)\ln(\mathcal{T}_h).
\label{eq.BLEND0}
\end{equation}
for $ h $ sufficiently small, 
thus expressing the derivative in terms of the logarithm of the shift operator.  

Expanding the (natural) logarithm in (\ref{eq.BLEND0}) as a formal (Taylor) power series around the identity operator $\mathcal{J}$ gives the following
formal representation:
\begin{equation}\label{op:eq}
\Delta=\frac{1}{h}\sum_{n=1}^\infty [(-1)^{n-1}(n-1)!] \frac{(\mathcal{T}_h-\mathcal{J})^n}{n!}=
-\frac{1}{h}\sum_{n=1}^\infty\frac{(\mathcal{J}-\mathcal{T}_h)^n}{n} ,
\end{equation}
for $ h $ sufficiently small.
Noting that the shift operator $\mathcal{T}_h$ satisfies 
$\mathcal{T}_h^k=\mathcal{T}_{kh},~k\geq 1$: 
\begin{equation}\label{eq:1}
(\mathcal{J}-\mathcal{T}_h)^n\phi=\sum_{k=0}^n \binom{n}{k} (-1)^k \mathcal{T}_{kh}\phi
=\sum_{k=0}^n(-1)^k\binom{n}{k}\phi(\cdot+kh),
\end{equation}
for $ n \geq 1 $.
Combining (\ref{op:eq}) with (\ref{eq:1}) yields
the {\em logarithmic series expansion of the derivative}:
\begin{equation}\label{eq:logser}
\Delta \phi ( \cdot ) =- \frac{1}{h}\sum_{n=1}^\infty 
\sum_{k=0}^n \frac{ (-1)^k  }{n} \binom{n}{k}\phi(\cdot+kh) ,
\end{equation}
for $  h$ sufficiently small.
A main contribution of this note is characterizing the domain of $ h $ for which 
(\ref{eq:logser}) will return the correct answer, 
a result established in Section~\ref{sec:BLEND}.

This logarithmic series expansion idea leading to (\ref{eq:logser})
was introduced by Asmussen and Glynn \cite[pp.212--213]{Glynn:07}
as a potential basis for FD estimators in the {\it Monte Carlo simulation} setting. 
For example, using just the first term in the series (\ref{eq:logser}) gives evaluation of two terms in (\ref{eq:1}) -- 
specifically $\Delta\approx(\mathcal{T}_h-\mathcal{J})/h$ in shorthand operator notation of (\ref{op:eq}), 
leading to a 
``simulation-based first-order approximation" 
that is just the usual forward difference estimator corresponding to (\ref{eq:forwardFD}). 
However, \cite[pp.213]{Glynn:07} concludes that due to numerical difficulties and variance challenges in the Monte Carlo simulation setting, 
% ``For all these reasons, 
``it is rare that estimators of higher order than the central difference estimator are used."

The main purpose of this note is to explore BLEND in the deterministic setting. 
Theoretical analysis and empirical evidence indicate that BLEND 
provides a practical alternative to other FD numerical approximations in many scenarios. 
In fact, the series expansion can also be used to derive {\it unbiased} simulation-based finite-difference estimators, but these suffer from numerical instabilities similar to those observed in  \cite{Glynn:07}.

The rest of this note is organized as follows. 
%The series expansion for the derivative is presented in Section~\ref{sec:RZ}.
The main theoretical results for BLEND are given in Section~\ref{sec:1a}, 
as well as the higher-dimensional extension to directional derivatives.
The BLEND algorithm is introduced in Section~\ref{sec:BLEND}.
Numerical examples are presented in Section~\ref{sec:2}, 
and Section~\ref{sec:stoch} concludes with a brief discussion of future application of BLEND to the stochastic setting.

\section{BLEND Analysis}\label{sec:1a}

The main technical conditions required for BLEND are the following: 

\begin{itemize}

\item[ {\bf (C1)} ]
The function $ \phi $  is analytical on $ ( \theta - \epsilon , \theta + h_0 + \epsilon ) $ for some $ \epsilon > 0 $.

\item[ {\bf (C2)} ]
Fix $ N$, and assume  that for $ n \leq N $:
$$
\sup_{ \theta \leq \hat \theta  \leq \theta + h n} \left|\phi^{(n)}(\hat \theta)\right|
 \leq  M b^n ,
$$
where $ \phi^{(0)} \equiv \phi $ and $ \phi^{(n)}$ denotes the $n$th derivative of $ \phi$.
\end{itemize}

\medskip

The condition on the analyticity of $ \phi $ on  $ ( \theta - \epsilon , \theta + h_0 + \epsilon ) $ ensures that all higher-order derivatives are well defined on the interval of interest $ [ \theta , \theta + h_0 ] $.
Note that if $ \phi $ is a polynomial of order $ k$, then {\bf (C2)} holds for $ b = \max ( \theta + h  k ,1 ) $ and $ M = k!$.
Conditions {\bf (C1)} and {\bf (C2)} allow a bounding of the $n^{th} $ power of the shift operator $ \left(\mathcal{J}-\mathcal{T}_h\right)^n$.

\begin{lemma}\label{le:bound}
If {\bf (C1)} and {\bf (C2)} hold, then for $ n \leq N $:
\[
\left | \left(\mathcal{J}-\mathcal{T}_h\right)^n\phi(\theta) 
\right | 
\leq     \frac{M }{\sqrt{2\pi n}}  \big (2 h \, b \,  e \big )^n .
\]
\end{lemma}
{\bf Proof:}
For any $n\geq N_h $ and $0\leq k\leq n$, we have a Taylor's series expansion
\begin{eqnarray*}
& & \phi(\theta+kh)  =\phi(\theta)+\frac{(kh)}{1!}\:\phi'(\theta)\\
& & \qquad \qquad +\frac{(kh)^2}{2!}\:\phi''(\theta)+\ldots+\frac{(kh)^n}{n!}\:\phi^{(n)}(\xi_k^n),
\end{eqnarray*}
for some $\xi_k^n\in(\theta,\theta+kh)$. Multiplying each expression by $(-1)^k\binom{n}{k}$ and summing up for
$k=1,\ldots,n$ yields (below, $\xi_k^i:=\theta$ for $0\leq i\leq n-1$ and any $k\leq n$):
\begin{eqnarray*}
\left(\mathcal{J}-\mathcal{T}_h\right)^n\phi(\theta) & = & \sum_{k=0}^n(-1)^k\binom{n}{k}\phi(\theta+kh) \\
& = & \sum_{k=0}^n(-1)^k\binom{n}{k}\sum_{i=0}^n\frac{(kh)^i}{i!}\phi^{(i)}(\xi_k^i) \\
& = & \sum_{i=0}^n\frac{h^i}{i!}\sum_{k=0}^n(-1)^k k^i \binom{n}{k} \phi^{(i)}(\xi_k^i) \\
& = & \frac{h^n}{n!}\sum_{k=0}^n(-1)^k k^n \binom{n}{k} \phi^{(n)}(\xi_k^n),
%=(-h)^n\phi^{(n)}(\xi_k^n),
\end{eqnarray*}
where, in the last equality, we have used the identity (Ruiz 1996, Corollary 2)
\begin{equation}\label{eq:haralambie}
\sum_{k=0}^n(-1)^k\binom{n}{k}\mathcal{P}(k)=0,
\end{equation}
%$$\sum_{k=0}^n(-1)^k\binom{n}{k}\mathcal{P}(k)=\sum_{k=0}^n(-1)^{n-k}\binom{n}{k}\mathcal{P}(n-k)=(-1)^n p_n n!,$$
valid for any polynomial $\mathcal{P}$ up to degree $n-1$. Hence, for $ n \geq 1 $:
\begin{eqnarray}\label{eq:bound}
\left|\left(\mathcal{J}-\mathcal{T}_h\right)^n\phi(\theta)\right|
 & \leq  & M\:\frac{h^n}{n!}\sum_{k=1}^n k^n \binom{n}{k}
\nonumber 
\\
& \leq & M\frac{(2n h)^n}{n!}\leq M\frac{(2h \, b \,  e)^n}{\sqrt{2\pi n}},
\end{eqnarray}
the last inequality following by Stirling's approximation. 
\hfill $ \Box $

\begin{example}\label{exp:ex}%\textrm{Exponential function}:
Let $\phi(\theta,x):=\theta\exp(-\theta x), x>0$. Then for any
$n\geq 1$ it holds that
\begin{eqnarray*}
& & \left(\mathcal{J}-\mathcal{T}_h\right)^n\phi(\theta,x) \\
&  & = \sum_{k=0}^n(-1)^k\binom{n}{k}(\theta+kh)\exp\left[-(\theta+kh)x\right] \\
&  & = \theta\exp(-\theta x)\sum_{k=0}^n(-1)^k\binom{n}{k}\exp(-khx) \\
& & \qquad \qquad 
+ h\exp(-\theta x)\sum_{k=1}^n(-1)^k k\binom{n}{k}\exp(-khx) \\
&  & = \theta\exp(-\theta x)\left[1-\exp(-hx)\right]^n  \\
& & \qquad - n h\exp\left[-(\theta+h)x\right]\left[1-\exp(-hx)\right]^{n-1}.
\end{eqnarray*}
Dividing by $-hn$ and summing everything up for $n\geq 1$, yields
\begin{equation}\label{erst}
-\frac{1}{h}\sum_{n\geq 1}\frac{\left(\mathcal{J}-\mathcal{T}_h\right)^n}{n}=\left(1-\theta x\right)\exp(-\theta x)=
\frac{\partial\phi}{\partial\theta}(\theta,x),
\end{equation}
which proves the validity of (\ref{op:eq}) for $\phi(\cdot,x)$.
\end{example}

Regarding the conditions of Lemma~\ref{le:bound}, we note that if $\phi$ is an analytic function defined only on a bounded interval $(l,b)$, then
one can use a change of variables $u:(a,\infty)\rightarrow (l,b)$ to obtain a new function $\varphi=\phi\circ u:(a,\infty)\rightarrow\mathbb{R}$ satisfying condition {\bf (C1)} and then recover the derivative $\phi'(\theta)$ by evaluating $\varphi'$ at the point $u^{-1}(\theta)$. Furthermore, recall that in general, analyticity amounts to the fact that for any compact set $D$, there exists $M,b>0$ such that
$$\sup_{\theta\in D}\left|\phi^{(n)}(\theta)\right|\leq M b^n n!.$$
 In that sense, condition {\bf (C2)} is more restrictive, as it imposes certain growth rates on increasing intervals.
Under the conditions of Lemma~\ref{le:bound}, we see that
$$\phi(\theta+h)=\sum_{n\geq 0}\frac{h^n}{n!}\:\phi^{(n)}(\theta)\leq M\exp(bh),$$
which suggests that the result applies to exponentially
bounded analytic functions.

%\subsection{Extension to Directional Derivatives}\label{sec:2b}

We conclude this section by discussing the extension to directional derivatives.
Let $ \theta \in \mathbb{R}^m $ and assume that $ \phi $ is analytical as a mapping of $ \theta $.
The directional derivative of $ \phi $ is direction $\vec{v}:=(v_1,\ldots,v_m) $ is given by 
\[
\nabla_{\vec{v}}\: \phi ( \theta ) = 
\lim_{ h \rightarrow 0 } \frac{   \phi ( \theta + h \vec{v})  - \phi ( \theta )}{ h || \vec{v} || }  = 
\sum_{i=1}^m  v_i \frac{\partial \phi}{\partial \theta_i } ( \theta ) ,
\]
where $ || \vec{v} || $ denotes the Euclidean norm.
Evaluating a directional derivative via the above analytical approach requires the evaluation of $m$ partial derivatives.
Alternatively, we can apply the shift-operator and obtain the following logarithmic series expansion of the directional derivative 
\begin{equation}\label{eq:logserdir}
\nabla_{\vec{v}}\: \phi ( \theta )  =- \frac{1}{h}\sum_{n=1}^\infty 
\sum_{k=0}^n \frac{ (-1)^k  }{n} \binom{n}{k}\phi(\theta +k  \vec{v} ) .
\end{equation}
Note that the complexity of the logarithmic series expansion is independent of the dimension $m$ of the parameter.

\section{The BLEND Algorithm}\label{sec:BLEND}

We consider the finite truncation of the sum in (\ref{op:eq})
as an approximate value for the derivative, i.e., 
the BLEND approximation of order $N$ and difference $h$ 
for the derivative of $ \phi ( \theta ) $ with respect to $ \theta $ is given by
\begin{equation}\label{eq:blend}
\Delta(N ,h )\phi ( \theta )  = - \frac{1}{h}\sum_{n=1}^N 
\sum_{k=0}^n \frac{ (-1)^k  }{n} \binom{n}{k}\phi(\theta +kh) ,
\end{equation}
and the error of the BLEND approximation will be denoted by
$$ 
R ( N , h, \phi ( \theta )) =  | \Delta  \phi ( \theta ) - \Delta (N,h) \phi ( \theta ) | .
$$

In the following we show how Lemma~\ref{le:bound} can be used to determine the truncation index $ N $ such that 
at least the first $ n$ digits of the derivative are exact.

\begin{lemma}\label{le:error}
Under conditions {\bf (C1)} and {\bf (C2)}, 
\[
R ( N , h ,\phi ( \theta )) \leq 
\frac{M }{\sqrt{2\pi } (N+1) ^{\frac{3}{2}} }   \frac{(2h \, b \, e)^{ (N+1)}}{1 - 2h b  e }, ~N \geq 1,
\]
provided $ h  < 1 / 2 b  e $.
\end{lemma}
{\bf Proof:}
 It follows from Lemma~\ref{le:bound} that 
\[
\forall n\geq 1:\:\left|\left(\mathcal{J}-\mathcal{T}_h\right)^n\phi(\theta)\right|
\leq M \frac{(2 hb_2 e)^n}{\sqrt{2\pi n}} .
\]
Hence, for $ h> 0 $ 
\begin{eqnarray*}
  \left | \sum_{n=N+1}^\infty \frac{(\mathcal{J}-\mathcal{T}_h)^n}{n} \right | 
& \leq & \frac{ M}{(N+1)^{\frac{3}{2}} \sqrt{2\pi}}  \sum_{n=N+1}^\infty  (2h b  e)^n  
\\
& = & \frac{ M}{(N+1)^{\frac{3}{2}} \sqrt{2\pi}}  (2hbe)^{N+1}   \frac{1}{1 - 2h b  e } ,
\end{eqnarray*}
provided $ h  < 1 / 2 b  e $. 
\hfill $ \Box$

\begin{example}
Consider  the mapping $ \phi ( \theta ) = \sin ( \theta )$.
We apply BLEND for computing the derivative of $ \sin ( \theta ) $ at $ \theta = 0 $. 
Condition {\bf (C2)} holds for all $ N$ with $ M =1 =b$.
By Lemma~\ref{le:error}, BLEND yields the correct result for $ h < 1 / 2 e \approx 0.1839 $.
\begin{table}[th!]
\begin{center}
\begin{tabular}{ l l l    }
N  & $ \Delta (N, 0.01) $  \\
\hline 
1   &             0.998334166468282             \\ 
2   &     {\bf 1.00}3321678961257   \\
3   &     {\bf   1.0000}29893016725  \\
4   &     {\bf     0.9999}80308400858 \\
5   &     {\bf    0.999999}646316608 \\
6   &     {\bf     1.000000}137620388  \\
7   &     {\bf    1.00000000}3815154  \\
8   &     {\bf    0.99999999}8963623 \\ \hline
 true & 1.0  \\ [8pt]
\end{tabular}
\caption{BLEND Approximation for the derivative $  \sin (\theta) $ at $ \theta=0$, small $h$}
\label{tab:approxsin}
\end{center}
\end{table}
As can be seen in Table~\ref{tab:approxsin} executing the BLEND algorithm for $ N=8 $ yields 8 exact digits of the derivative.
To illustrate the impact on the bound on $ h$ put forward in Lemma~\ref{le:error}, we present in 
Table~\ref{tab:approxsinb} the output of BLEND for $ h =1$. 

\begin{table}[th!]
\begin{center}
\begin{tabular}{ l l l    }
N  & $ \Delta (N, 1.0) $  \\
\hline 
1   &     0.841470984807897    \\ 
2   &      1.228293256202952  \\
3   &     1.207506816871789   \\
4   &       1.015352293328013 \\
5   &      0.885486080979581 \\
6   &      0.903764738896000   \\
7   &    1.003453862663737  \\
8   &      1.071046882890327 \\ \hline
\vspace*{0.1pt}
\end{tabular}
\caption{BLEND Approximation for the derivative $  \sin (\theta) $ at $ \theta=0$, large $h$}
\label{tab:approxsinb}
\end{center}
\end{table}

As can be seen in Table~\ref{tab:approxsinb} executing the BLEND algorithm for $ h=1 $ fails to produce the correct output, which stems from the fact that $ 1 > 1 /2 e \approx 0.1839 $.

The fact that BLEND fails to yield the correct output for large $ h$  can best be seen by applying BLEND to the derivative of $ \sin (\theta ) $ at $ \theta = 0 $ with $ h = 2\pi $, which yields $ 0   = \Delta \phi (N , 2 \pi  ) $ for all $ N$, whereas  $ 1 $ is the correct answer. 
\end{example}

Lemma~\ref{le:error} shows that the error of considering $ \Delta ( N , h  ) \phi ( \theta ) $ rather than $ \Delta (h) \phi ( \theta )  $ 
is of order $ O ((2 h b  e )^{(N+1)} )$.
Moreover, for given $ N$, provided $ h $ is sufficiently small so that it satisfies the condition in Lemma~\ref{le:error} and upper bounds for $ M $ and $ b$ are known, solving for $ h $ in 
\begin{equation}\label{eq:exactfd}
\frac{M }{\sqrt{2\pi } 2^{\frac{N+1}{2}} }   \frac{(2h \, b \, e)^{N+1}}{1 - 2h b  e } = 10^{-(K+1)}
\end{equation}
yields a value for $ h $ such that the FWD (forward finite difference) approximation  $ \Delta ( N , h  ) \phi ( \theta )$  is exact in at least the first $ K$ digits. 
We call this {\em the $ K$-exact FWD approximation}.

Typically, $ M $ and $b$ cannot be computed exactly. 
In this case the bound on the remainder put forward in Lemma~\ref{le:error}  can be facilitated through the geometric convergence rate of the series
in (\ref{eq:blend}) for small $ h$.
Specifically, inspecting the values of $ \Delta ( N , h ) $ for $ h $ fixed and $ N = 1 , 2, \ldots $ one expects that after a possible transient phase the geometric error decrease established in  Lemma~\ref{le:error} sets in.
This gives rise to the following heuristic: Compute  $ \Delta ( N , h ) $, for $ N = 1 , 2 , \ldots, N_{\max} $, with $ N_{\max} $ a pre-specified number. 
If no stabilization is seen, then decrease $h$ and compute again the $N_{\max}$ values. If stabilized with the first $L$ digits being the same for $N_{\max-1}$ and $N_{\max}$, then
in light of the geometric error decrease, accept the first $ L $ digits in   $ \Delta ( N' , h ) $ -- as if the values of  $ \Delta ( N , h ) $ for $ N > N'$ only affected digits $ L+1 $ and larger.
We call this the {\em stabilization indication} of the series.
We summarize this in BLEND algorithm presented in the following.\\

\noindent
{\bf BLEND Algorithm}
{\em \begin{itemize}
\item[(i)] 
Choose $ h $ small and set $ N_{\max} $.

\item[(ii)] Evaluate $ \Delta(N ,h )\phi ( \theta )  $ in (\ref{eq:blend}) for $ N =1 , \ldots , N_{\max} $.

\item[(iii)]  Let $ L $ be such that the first $ L $ digits in $ \Delta (N_{\max} -1 ,h )\phi ( \theta )  $ and  $ \Delta( N_{\max} ,h )\phi ( \theta )  $ are equal, then 
the value of $ \Delta(N_{\max} ,h )\phi ( \theta )  $ up to  first $ L $ digits is outputted as the exact value for the first $ L $ digits of the derivative.

\end{itemize}
}

\begin{remark}
The BLEND exploits parallel processing by executing each of the evaluations in step (ii) simultaneously. 
As an illustration, $N_{\max}=8 $ would be a natural choice for many current portable multi-core computing platforms that are commonly available, e.g., Apple laptops.   
\end{remark}

\section{Numerical Examples}\label{sec:2}

We illustrate the $K$-exact FD approximation as well as the BLEND algorithm for a series of examples, 
beginning with some toy examples and then considering a more complicated example that illustrates the type of setting 
that we envision for the main application of BLEND.  
For comparisons, recall that the FD approximation is just the first term in the tables below, when $N=1$.

\begin{example}
Let $ \phi ( \theta ) = 5 \theta^4 $ and note that {\bf (C1)} holds for $ h_0 = \infty $. 
We will compute the derivative at $ \theta_0 = 2$.
In order to apply condition {\bf (C2)}, we assume that  $ h $ is no larger than $ \hat h = 0.1 $.
Only the first 4 derivatives are significant and we arrive at
\[
\sup_{ \theta_0 \leq  \theta  \leq \theta_0 + \hat h n} \left|\phi^{(n)}(\theta)\right| 
 \leq 120 (2 + 0.1 \times 4 )^4,  1 \leq n \leq 4 ,
 \]
and thus $ M =120 $ and  $ b = 2.4 $.

We apply BLEND for $ N_{\max } =8 $ and $ h = 0.01 $.
The numerical results are provided in Table~\ref{tab:approx3}.
\begin{table}[th!]
\begin{center}
\begin{tabular}{ l l l    }
N  & $ \Delta (N, 0.01) $  \\
\hline 
1   &      160.1200400049834       \\ 
2   &      {\bf 159.9999}199699909    \\
3   &      {\bf 160.000}0299999870 \\
4   &      {\bf 159.9999999999}799  \\
5    &     {\bf  159.99999999997}19  \\
6   &      {\bf  159.9999999999}577   \\
	7   &      {\bf  159.9999999999}281\\
8   &      {\bf  159.9999999999}981   \\ \hline
 true &  160   \\ [8pt]
\end{tabular}
\caption{The Finite Approximation for the derivative $ 5 \theta^4 $ at $ \theta=2$, $h=0.01$}
\label{tab:approx3}
\end{center}
\end{table}
The numerical results put forward in Table~\ref{tab:approx3} show that our stabilization indication works in natural way.
Notice that up to 10 digits  $  \Delta ( 8, 0.01) $  is rounded off to $160$, giving same digits up to precision $10^{-10}$. 

We now turn to $ K$-exact FWD approximation.
Suppose we want to compute the first 6 digits of the derivative exactly, with $N_{\max}=2$.
Then we solve (\ref{eq:exactfd}) for $ h$ with $ M =120 $ and  $ b = 2.4 $, i.e., 
\[
\frac{120 }{\sqrt{2\pi } 2^{\frac{3}{2}} }   \frac{(4.8 h \, e)^{ 3}}{1 - 4.8  h  e } = 10^{-7 },
\]
so that  $ h $ has to be smaller than $ 0.00013$. 
The resulting value for $ \Delta ( 2 , 0.00013)$ is $ 1.599999999999915$.
\end{example}

In the following we illustrate the application of our results to directional derivatives.

\begin{example}
Consider $ \theta \in \mathbb{R}^m$ and let
\[
\phi ( \theta ) = \sum_{ i=1}^m  a_i \theta_i^2   ,
\]
for some constants $ a_i \in \mathbb{R}$, $ 1 \leq i \leq m $.
We consider the directional derivative of $ \phi ( \theta ) $ in direction $ \vec{v} $, where
$ \vec{v} \in [-1 , 1 ]^m $ is chosen such that $ \sum_{i=1}^m v_i = 0 $.

For our numerical experiment, we let $ m = 9 $, $ a_i = 2^{-i}$, for $ 1 \leq i \leq m $,
and 
$$ 
\vec{v}=  (1/3) ( -1 , 1,  -1, -1 , 1 , 1 , 1 , -1 , 1 ), 
$$ so that $ || \vec {v} || =1 $.
We obtain
\[
\frac{\partial\phi}{\partial \theta_i} ( \theta ) = \frac{\theta_i}{2^{-i+1}}  , \quad \hbox{ for }  1 \leq i \leq m .
\]
We now apply BLEND for $ N_{\max } =8 $ and $ h = 0.01 $ for computing the 
directional derivative of $ \phi ( \theta) $ in direction $ \vec{v} $ at 
$ \theta $, with $ \theta_i = i$, for $ 1 \leq i \leq m $.
The numerical results are provided in Table~\ref{tab:direct}.
\begin{table}[th!]
\begin{center}
\begin{tabular}{ l l l    }
N  & $ \Delta (N, 0.001) $  \\
\hline 
1   &           3.958029296875054      \\ 
2   &     {\bf 3.95}7031250000576      \\
3   &     {\bf 3.95703125000}2056 \\
4   &     {\bf 3.95703125000}4276       \\
5    &    {\bf 3.95703125000}5520  \\
6   &     {\bf  3.95703125000}7444  \\
7   &     {\bf  3.9570312500}13154  \\
8   &     {\bf  3.957031250013}043 \\ \hline
 true &  3.9570312500138101  \\ [8pt]
\end{tabular}
\caption{The BLEND Finite Approximation for the directional derivative}
\label{tab:direct}
\end{center}
\end{table}

\end{example}

Now we consider the queueing example,
where the constants required to set the $K$-exact FD approximation would not be available.

\begin{example}

%\subsection{Stationary Markov Chains}\label{ex:MC}

Consider a two-station tandem queueing system with finite capacity $ N_1 $ at station 1 and $ N_2 $ at station 2.
Jobs arrive to the network according to a Poisson process with arrival rate $ \lambda$, and service times are 
independent and identically distributed (i.i.d.) exponential with rate $ \mu_i $  at station $i=1,2$.
When there is no waiting place available at station 1, an arrival is rejected and lost.
When there is no waiting place available at station 2, service station 1 is stopped.
This example is taken from \cite{DDD}, and we refer for motivation and details to the references therein.
Due to the finite buffers, no closed-form expression for the stationary distribution exists.
However, letting $ Q $ denote the infinitesimal generator of the process, the stationary distribution solves $ \pi Q = 0$ with normalizing equation $ \sum \pi_i =1 $ and is easily numerically evaluated.

We apply the  BLEND algorithm for computing the derivative of the blocking probability with respect to  $ \lambda$,
taking $ \lambda =1 $, $ \mu=1$, $ \eta_1 = 1 $ and $ \eta_2 =2 $, 
with finite capacity queue sizes $ N_1= N_2 =10$.
The numerical results are given in Table~\ref{tab:tandem}, 
where the true value has been obtained by the finite difference method through a series of experiments.

\begin{table}[th!]
\begin{center}
\begin{tabular}{ l l l    }
N  & $ \Delta (N, 0.01 ) $  \\
\hline 
1   &                 0.613180514116096       \\ 
2   &       {\bf    0.61}0046682208255  \\
3   &      {\bf     0.609}671969013386  \\
4   &      {\bf     0.6096}61671019043 \\
5   &      {\bf    0.60966}2935724646 \\
6  &        {\bf    0.60966}3162694883  \\
7   &     {\bf    0.6096631}73459084 \\
8  &    {\bf     0.60966317}0509458 \\ \hline
 true &  0.609663168  \\ [8pt]
\end{tabular}
\caption{The BLEND  Approximation for the Loss Probability Sensitivity with respect to $ \lambda $ at $ \lambda =1 $}
\label{tab:tandem}
\end{center}
\end{table}

\end{example}

\section{Conclusion and Future Research} 
\label{sec:stoch}

We presented a new finite difference derivative approximation called the BLEND algorithm.
BLEND is particularly useful when the expression of $\phi$ (the function of interest) is not available in closed-form, but the values at arbitrary points within its domain can be efficiently numerically computed. 
We characterized the value of the difference parameter $h$ for which the BLEND algorithm applies, 
and provided a bound on the error that leads to an estimate on the number of terms required to achieve a particular degree of precision. 
Various numerical examples illustrate the effectiveness of the BLEND algorithm, including the importance of choosing the parameters correctly and the practical implementation using parallel computing. 

Future research of interest is considering the extension to the stochastic case by letting $ \phi ( \theta ) = \mathbb{E} [ Z ( \theta ) ] $, for $ Z ( \theta ) $ the underlying stochastic variable, and the application of BLEND to high dimensional problems.

\end{document}